\documentclass{amsart}
\usepackage{amsfonts}

\newtheorem{theorem}{Theorem}
\theoremstyle{plain}

\newtheorem{corollary}{Corollary}

\newtheorem{lemma}{Lemma}

\newtheorem{remark}{Remark}

\numberwithin{equation}{section}
\input{tcilatex}

\begin{document}
\title[Bombieri Inequality]{Some Bombieri Type Inequalities in Inner Product Spaces}
\author{S.S. Dragomir}
\address{School of Computer Science and Mathematics\\
Victoria University of Technology\\
PO Box 14428, MCMC \\
Victoria 8001, Australia.}
\email{sever.dragomir@vu.edu.au}
\urladdr{http://rgmia.vu.edu.au/SSDragomirWeb.html}
\date{11 June, 2003.}

\begin{abstract}
Companion results to the Bombieri generalisation of Bessel's inequality in
inner product spaces are given.
\end{abstract}

\maketitle

\section{Introduction}

In 1971, E. Bombieri \cite{1b}, has given the following generalisation of
Bessel's inequality: 
\begin{equation}
\sum_{i=1}^{n}\left| \left( x,y_{i}\right) \right| ^{2}\leq \left\|
x\right\| ^{2}\max_{1\leq i\leq n}\left\{ \sum_{j=1}^{n}\left| \left(
y_{i},y_{j}\right) \right| \right\} ,  \label{1.1}
\end{equation}
where $x,y_{1},\dots ,y_{n}$ are vectors in the inner product space $\left(
H;\left( \cdot ,\cdot \right) \right) .$

It is obvious that if $\left( y_{i}\right) _{1\leq i\leq n}=\left(
e_{i}\right) _{1\leq i\leq n},$ where $\left( e_{i}\right) _{1\leq i\leq n}$
are orthornormal vectors in $H,$ i.e., $\left( e_{i},e_{j}\right) =\delta
_{ij}$ $\left( i,j=1,\dots ,n\right) ,$ where $\delta _{ij}$ is the
Kronecker delta, then (\ref{1.1}) provides Bessel's inequality 
\begin{equation}
\sum_{i=1}^{n}\left| \left( x,e_{i}\right) \right| ^{2}\leq \left\|
x\right\| ^{2},\ \ \ x\in H.  \label{1.2}
\end{equation}

In this paper we point out other Bombieri type inequalities and show that,
some times, the new ones may provide better bounds for $\sum_{i=1}^{n}\left%
\vert \left( x,y_{i}\right) \right\vert ^{2}.$

\section{The Results}

The following lemma which is of interest in itself holds.

\begin{lemma}
\label{l2.1}Let $z_{1},\dots ,z_{n}\in H$ and $\alpha _{1},\dots ,\alpha
_{n}\in \mathbb{K}$. Then one has the inequalities: 
\begin{eqnarray}
&&\left\| \sum_{i=1}^{n}\alpha _{i}z_{i}\right\| ^{2}  \label{2.1a} \\
&\leq &\left\{ 
\begin{array}{l}
\max\limits_{1\leq i\leq n}\left| \alpha _{i}\right|
^{2}\sum_{i,j=1}^{n}\left| \left( z_{i},z_{j}\right) \right| ; \\ 
\\ 
\left( \sum_{i=1}^{n}\left| \alpha _{i}\right| ^{p}\right) ^{\frac{2}{p}%
}\left( \sum_{i,j=1}^{n}\left| \left( z_{i},z_{j}\right) \right| ^{q}\right)
^{\frac{1}{q}},\ \ \text{where \ }p>1,\ \frac{1}{p}+\frac{1}{q}=1; \\ 
\\ 
\left( \sum_{i=1}^{n}\left| \alpha _{i}\right| \right) ^{2}\
\max\limits_{1\leq i,j\leq n}\left| \left( z_{i},z_{j}\right) \right| ;
\end{array}
\right.  \notag \\
&\leq &\left\{ 
\begin{array}{l}
\max\limits_{1\leq i\leq n}\left| \alpha _{i}\right| ^{2}\left(
\sum_{i=1}^{n}\left\| z_{i}\right\| \right) ^{2}; \\ 
\\ 
\left( \sum_{i=1}^{n}\left| \alpha _{i}\right| ^{p}\right) ^{\frac{2}{p}%
}\left( \sum_{i=1}^{n}\left\| z_{i}\right\| ^{q}\right) ^{\frac{2}{q}},\ \ 
\text{where \ }p>1,\ \frac{1}{p}+\frac{1}{q}=1; \\ 
\\ 
\left( \sum_{i=1}^{n}\left| \alpha _{i}\right| \right) ^{2}\
\max\limits_{1\leq i\leq n}\left\| z_{i}\right\| ^{2}.
\end{array}
\right.  \notag
\end{eqnarray}
\end{lemma}

\begin{proof}
We observe that 
\begin{align}
\left\| \sum_{i=1}^{n}\alpha _{i}z_{i}\right\| ^{2}& =\left(
\sum_{i=1}^{n}\alpha _{i}z_{i},\sum_{j=1}^{n}\alpha _{j}z_{j}\right)
=\sum_{i=1}^{n}\sum_{j=1}^{n}\alpha _{i}\overline{\alpha _{j}}\left(
z_{i},z_{j}\right)  \label{2.2} \\
& =\left| \sum_{i=1}^{n}\sum_{j=1}^{n}\alpha _{i}\overline{\alpha _{j}}%
\left| \left( z_{i},z_{j}\right) \right| \right| \leq
\sum_{i=1}^{n}\sum_{j=1}^{n}\left| \alpha _{i}\right| \left| \alpha
_{j}\right| \left| \left( z_{i},z_{j}\right) \right| =:M.  \notag
\end{align}
Firstly, we have 
\begin{align*}
M& \leq \max\limits_{1\leq i,j\leq n}\left\{ \left| \alpha _{i}\right|
\left| \alpha _{j}\right| \right\} \sum\limits_{i,j=1}^{n}\left| \left(
z_{i},z_{j}\right) \right| \\
& =\max\limits_{1\leq i\leq n}\left| \alpha _{i}\right|
^{2}\sum\limits_{i,j=1}^{n}\left| \left( z_{i},z_{j}\right) \right| .
\end{align*}
Secondly, by the H\"{o}lder inequality for double sums, we have 
\begin{align*}
M& \leq \left[ \sum\limits_{i,j=1}^{n}\left( \left| \alpha _{i}\right|
\left| \alpha _{j}\right| \right) ^{p}\right] ^{\frac{1}{p}}\left(
\sum\limits_{i,j=1}^{n}\left| \left( z_{i},z_{j}\right) \right| ^{q}\right)
^{\frac{1}{q}} \\
& =\left( \sum_{i=1}^{n}\left| \alpha _{i}\right| ^{p}\sum_{j=1}^{n}\left|
\alpha _{j}\right| ^{p}\right) ^{\frac{1}{p}}\left(
\sum\limits_{i,j=1}^{n}\left| \left( z_{i},z_{j}\right) \right| ^{q}\right)
^{\frac{1}{q}} \\
& =\left( \sum\limits_{i=1}^{n}\left| \alpha _{i}\right| ^{p}\right) ^{\frac{%
2}{p}}\left( \sum\limits_{i,j=1}^{n}\left| \left( z_{i},z_{j}\right) \right|
^{q}\right) ^{\frac{1}{q}},
\end{align*}
where $p>1,$ $\frac{1}{p}+\frac{1}{q}=1.$

Finally, we have 
\begin{equation*}
M\leq \max\limits_{1\leq i,j\leq n}\left| \left( z_{i},z_{j}\right) \right|
\sum\limits_{i,j=1}^{n}\left| \alpha _{i}\right| \left| \alpha _{j}\right|
=\left( \sum\limits_{i=1}^{n}\left| \alpha _{i}\right| \right)
^{2}\max\limits_{1\leq i,j\leq n}\left| \left( z_{i},z_{j}\right) \right|
\end{equation*}
and the first part of the lemma is proved.

The second part is obvious on taking into account, by Schwarz's inequality
in $H$, that we have 
\begin{equation*}
\left| \left( z_{i},z_{j}\right) \right| \leq \left\| z_{i}\right\| \left\|
z_{j}\right\| ,
\end{equation*}
for any $i,j\in \left\{ 1,...,n\right\} .$ We omit the details.
\end{proof}

\begin{corollary}
\label{c2.2}With the assumptions in Lemma \ref{l2.1}, one has 
\begin{eqnarray}
\left\| \sum_{i=1}^{n}\alpha _{i}z_{i}\right\| ^{2} &\leq
&\sum\limits_{i=1}^{n}\left| \alpha _{i}\right| ^{2}\left(
\sum\limits_{i,j=1}^{n}\left| \left( z_{i},z_{j}\right) \right| ^{2}\right)
^{\frac{1}{2}}  \label{2.3} \\
&\leq &\sum\limits_{i=1}^{n}\left| \alpha _{i}\right|
^{2}\sum_{i=1}^{n}\left\| z_{i}\right\| ^{2}.  \notag
\end{eqnarray}
\end{corollary}

The proof follows by Lemma \ref{l2.1} on choosing $p=q=2.$

Note also that $\left( \ref{2.3}\right) $ provides a refinement of the well
known Cauchy-Bunyakovsky-Schwarz inequality for sequences of vectors in
inner product spaces, namely 
\begin{equation*}
\left\| \sum_{i=1}^{n}\alpha _{i}z_{i}\right\| ^{2}\leq
\sum\limits_{i=1}^{n}\left| \alpha _{i}\right| ^{2}\sum_{i=1}^{n}\left\|
z_{i}\right\| ^{2}.
\end{equation*}

The following lemma also holds.

\begin{lemma}
\label{l2.4}Let $x,y_{1},\dots ,y_{n}\in H$ and $c_{1},\dots ,c_{n}\in 
\mathbb{K}$. Then one has the inequalities: 
\begin{eqnarray}
&&\left| \sum\limits_{i=1}^{n}c_{i}\left( x,y_{i}\right) \right| ^{2}
\label{2.4a} \\
&\leq &\left\| x\right\| ^{2}\times \left\{ 
\begin{array}{l}
\ \max\limits_{1\leq i\leq n}\left| c_{i}\right|
^{2}\sum\limits_{i,j=1}^{n}\left| \left( y_{i},y_{j}\right) \right| ; \\ 
\\ 
\left( \sum\limits_{i=1}^{n}\left| c_{i}\right| ^{p}\right) ^{\frac{2}{p}%
}\left( \sum\limits_{i,j=1}^{n}\left| \left( y_{i},y_{j}\right) \right|
^{q}\right) ^{\frac{1}{q}},\ \ \ \ \text{where}\ \ p>1,\ \frac{1}{p}+\frac{1%
}{q}=1; \\ 
\\ 
\left( \sum\limits_{i=1}^{n}\left| c_{i}\right| \right)
^{2}\max\limits_{1\leq i,j\leq n}\left| \left( y_{i},y_{j}\right) \right| ;
\end{array}
\right.  \notag \\
&\leq &\left\| x\right\| ^{2}\times \left\{ 
\begin{array}{l}
\ \max\limits_{1\leq i\leq n}\left| c_{i}\right| ^{2}\left(
\sum\limits_{i=1}^{n}\left\| y_{i}\right\| \right) ^{2}; \\ 
\\ 
\left( \sum\limits_{i=1}^{n}\left| c_{i}\right| ^{p}\right) ^{\frac{2}{p}%
}\left( \sum\limits_{i=1}^{n}\left\| y_{i}\right\| ^{q}\right) ^{\frac{2}{q}%
},\ \ \ \ \text{where}\ \ p>1,\ \frac{1}{p}+\frac{1}{q}=1; \\ 
\\ 
\left( \sum\limits_{i=1}^{n}\left| c_{i}\right| \right)
^{2}\max\limits_{1\leq i\leq n}\left\| y_{i}\right\| ^{2}.
\end{array}
\right.  \notag
\end{eqnarray}
\end{lemma}

\begin{proof}
We have, by Schwarz's inequality in the inner product $\left( H;\left( \cdot
,\cdot \right) \right) ,$ that 
\begin{equation*}
\left| \sum\limits_{i=1}^{n}c_{i}\left( x,y_{i}\right) \right| ^{2}=\left|
\left( x,\sum\limits_{i=1}^{n}\overline{c_{i}}y_{i}\right) \right| ^{2}\leq
\left\| x\right\| ^{2}\left\| \sum\limits_{i=1}^{n}\overline{c_{i}}%
y_{i}\right\| ^{2}.
\end{equation*}
Now, applying Lemma \ref{l2.1} for $\alpha _{i}=\overline{c_{i}}$, $%
z_{i}=y_{i}$ $\left( i=1,\dots ,n\right) ,$ the inequality (\ref{2.4a}) is
proved.
\end{proof}

\begin{corollary}
\label{c2.5}With the assumptions in Lemma \ref{l2.4}, one has 
\begin{eqnarray}
\left| \sum\limits_{i=1}^{n}c_{i}\left( x,y_{i}\right) \right| ^{2} &\leq
&\left\| x\right\| ^{2}\sum\limits_{i=1}^{n}\left| c_{i}\right| ^{2}\left(
\sum\limits_{i,j=1}^{n}\left| \left( y_{i},y_{j}\right) \right| ^{2}\right)
^{\frac{1}{2}}  \label{2.5} \\
&\leq &\left\| x\right\| ^{2}\sum\limits_{i=1}^{n}\left| c_{i}\right|
^{2}\sum\limits_{i=1}^{n}\left\| y_{i}\right\| ^{2}.  \notag
\end{eqnarray}
\end{corollary}

The proof follows by Lemma \ref{l2.4}, on choosing $p=q=2.$

\begin{remark}
\label{r2.6}The inequality (\ref{2.5}) was firstly obtained in \cite{2b}
(see inequality (7)).
\end{remark}

The following theorem incorporating three Bombieri type inequalities holds.

\begin{theorem}
\label{t2.7}Let $x,y_{1},\dots ,y_{n}\in H.$ Then one has the inequalities: 
\begin{multline}
\sum\limits_{i=1}^{n}\left| \left( x,y_{i}\right) \right| ^{2}  \label{2.6}
\\
\leq \left\| x\right\| \times \left\{ 
\begin{array}{l}
\ \max\limits_{1\leq i\leq n}\left| \left( x,y_{i}\right) \right| \left(
\sum\limits_{i,j=1}^{n}\left| \left( y_{i},y_{j}\right) \right| \right) ^{%
\frac{1}{2}}; \\ 
\\ 
\left( \sum\limits_{i=1}^{n}\left| \left( x,y_{i}\right) \right| ^{p}\right)
^{\frac{1}{p}}\left( \sum\limits_{i,j=1}^{n}\left| \left( y_{i},y_{j}\right)
\right| ^{q}\right) ^{\frac{1}{2q}},\ \ \ \ \text{where}\ \ p>1,\ \frac{1}{p}%
+\frac{1}{q}=1; \\ 
\\ 
\sum\limits_{i=1}^{n}\left| \left( x,y_{i}\right) \right| \
\max\limits_{1\leq i,j\leq n}\left| \left( y_{i},y_{j}\right) \right| ^{%
\frac{1}{2}}.
\end{array}
\right.
\end{multline}
\end{theorem}

\begin{proof}
Choosing $c_{i}=\overline{\left( x,y_{i}\right) }$ $\left( i=1,\dots
,n\right) $ in (\ref{2.4a}) we deduce 
\begin{multline}
\left( \sum\limits_{i=1}^{n}\left| \left( x,y_{i}\right) \right| ^{2}\right)
^{2}  \label{2.7} \\
\leq \left\| x\right\| ^{2}\times \left\{ 
\begin{array}{l}
\ \max\limits_{1\leq i\leq n}\left| \left( x,y_{i}\right) \right| ^{2}\left(
\sum\limits_{i,j=1}^{n}\left| \left( y_{i},y_{j}\right) \right| \right) ; \\ 
\\ 
\left( \sum\limits_{i=1}^{n}\left| \left( x,y_{i}\right) \right| ^{p}\right)
^{\frac{2}{p}}\left( \sum\limits_{i,j=1}^{n}\left| \left( y_{i},y_{j}\right)
\right| ^{q}\right) ^{\frac{1}{q}},\ \ \ \ \text{where}\ \ p>1,\ \frac{1}{p}+%
\frac{1}{q}=1; \\ 
\\ 
\left( \sum\limits_{i=1}^{n}\left| \left( x,y_{i}\right) \right| \right)
^{2}\ \max\limits_{1\leq i,j\leq n}\left| \left( y_{i},y_{j}\right) \right| ;
\end{array}
\right.
\end{multline}
which, by taking the square root, is clearly equivalent to (\ref{2.6}).
\end{proof}

\begin{remark}
If $\left( y_{i}\right) _{1\leq i\leq n}=\left( e_{i}\right) _{1\leq i\leq
n} $ where $\left( e_{i}\right) _{1\leq i\leq n}$ are orthornormal vectors
in $H,$ then by (\ref{2.6}) we deduce 
\begin{equation}
\sum\limits_{i=1}^{n}\left| \left( x,e_{i}\right) \right| ^{2}\leq \left\|
x\right\| \times \left\{ 
\begin{array}{l}
\sqrt{n}\ \max\limits_{1\leq i\leq n}\left| \left( x,e_{i}\right) \right| ;
\\ 
\\ 
n^{\frac{1}{2q}}\left( \sum\limits_{i=1}^{n}\left| \left( x,e_{i}\right)
\right| ^{p}\right) ^{\frac{1}{p}},\ \ \ \ \text{where}\ \ p>1,\ \frac{1}{p}+%
\frac{1}{q}=1; \\ 
\\ 
\sum\limits_{i=1}^{n}\left| \left( x,e_{i}\right) \right| .
\end{array}
\right.  \label{2.7.a}
\end{equation}
\end{remark}

If in (\ref{2.7}) we take $p=q=2,$ then we obtain the following inequality
that was formulated in \cite[p. 81]{2b}.

\begin{corollary}
\label{c2.8}With the assumptions in Theorem \ref{t2.7}, we have: 
\begin{equation}
\sum\limits_{i=1}^{n}\left\vert \left( x,y_{i}\right) \right\vert ^{2}\leq
\left\Vert x\right\Vert ^{2}\left( \sum\limits_{i,j=1}^{n}\left\vert \left(
y_{i},y_{j}\right) \right\vert ^{2}\right) ^{\frac{1}{2}}.  \label{2.8}
\end{equation}
\end{corollary}

\begin{remark}
\label{r2.9}Observe, that by the monotonicity of power means, we may write 
\begin{equation}
\left( \frac{\sum_{i=1}^{n}\left| \left( x,y_{i}\right) \right| ^{p}}{n}%
\right) ^{\frac{1}{p}}\leq \left( \frac{\sum_{i=1}^{n}\left| \left(
x,y_{i}\right) \right| ^{2}}{n}\right) ^{\frac{1}{2}},\ \ 1<p\leq 2.
\end{equation}
Taking the square in both sides, one has 
\begin{equation*}
\left( \frac{\sum_{i=1}^{n}\left| \left( x,y_{i}\right) \right| ^{p}}{n}%
\right) ^{\frac{2}{p}}\leq \frac{\sum_{i=1}^{n}\left| \left( x,y_{i}\right)
\right| ^{2}}{n},
\end{equation*}
giving 
\begin{equation}
\left( \sum_{i=1}^{n}\left| \left( x,y_{i}\right) \right| ^{p}\right) ^{%
\frac{2}{p}}\leq n^{\frac{2}{p}-1}\sum_{i=1}^{n}\left| \left( x,y_{i}\right)
\right| ^{2}.  \label{2.10}
\end{equation}
Using (\ref{2.10}) and the second inequality in (\ref{2.7}) we may deduce
the following result 
\begin{equation}
\sum_{i=1}^{n}\left| \left( x,y_{i}\right) \right| ^{2}\leq n^{\frac{2}{p}%
-1}\left\| x\right\| ^{2}\left( \sum\limits_{i,j=1}^{n}\left| \left(
y_{i},y_{j}\right) \right| ^{q}\right) ^{\frac{1}{q}},  \label{2.11}
\end{equation}
for $1<p\leq 2,$ $\frac{1}{p}+\frac{1}{q}=1.$

Note that for $p=2$ $\left( q=2\right) $ we recapture (\ref{2.8}).
\end{remark}

\begin{remark}
Let us compare Bombieri's result 
\begin{equation}
\sum\limits_{i=1}^{n}\left| \left( x,y_{i}\right) \right| ^{2}\leq \left\|
x\right\| ^{2}\max\limits_{1\leq i\leq n}\left\{ \sum\limits_{j=1}^{n}\left|
\left( y_{i},y_{j}\right) \right| \right\}  \label{2.12}
\end{equation}
with our general result (\ref{2.11}).

To do that, denote 
\begin{equation*}
M_{1}:=\max\limits_{1\leq i\leq n}\left\{ \sum\limits_{j=1}^{n}\left| \left(
y_{i},y_{j}\right) \right| \right\}
\end{equation*}
and 
\begin{equation*}
M_{2}:=n^{\frac{2}{p}-1}\left( \sum\limits_{i,j=1}^{n}\left| \left(
y_{i},y_{j}\right) \right| ^{q}\right) ^{\frac{1}{q}},\ \ \ 1<p\leq 2,\ 
\frac{1}{p}+\frac{1}{q}=1.
\end{equation*}
Consider the inner product space $H=\mathbb{R}$, $\left( x,y\right) =x\cdot
y $, $n=2$ and $y_{1}=a>0,$ $y_{2}=b>0.$ Then 
\begin{gather*}
M_{1}=\max \left\{ a^{2}+ab,ab+b^{2}\right\} =\left( a+b\right) \max \left\{
a,b\right\} , \\
M_{2}=2^{\frac{2}{p}-1}\left( a^{q}+b^{q}\right) ^{\frac{2}{q}}=2^{\frac{2}{p%
}-1}\left( a^{\frac{p}{p-1}}+b^{\frac{p}{p-1}}\right) ^{\frac{2\left(
p-1\right) }{p}},\ 1<p\leq 2.
\end{gather*}
Assume $a=1,$ $b\in \left[ 0,1\right] ,$ $p\in (1,2].$ Utilizing Maple 6,
one may easily see by plotting the function 
\begin{equation*}
f\left( b,p\right) :=M_{2}-M_{1}=2^{\frac{2}{p}-1}\left( 1+b^{\frac{p}{p-1}%
}\right) ^{\frac{2\left( p-1\right) }{p}}-1-b
\end{equation*}
that it has positive and negative values in the box $\left[ 0,1\right]
\times \left[ 1,2\right] $, showing that the inequalities (\ref{2.11}) and (%
\ref{2.12}) cannot be compared. This means that one is not always better
than the other.
\end{remark}

\end{document}